\documentclass[a4paper,leqno,12pt]{article}
\usepackage{amsmath,amsfonts,amssymb}
\usepackage[dvips,%
pdftitle={A Translation of the T. Levi-Civita paper:
Interpretazione Gruppale
    degli Integrali di un Sistema Canonico},%
pdfauthor={G. Saccomandi, R. Vitolo},%
pdfkeywords={Levi-Civita, point symmetries, conservation laws,
  Hamilton's equations, Noether's theorem.}%
]{hyperref}

\newcommand{\be}{\begin{equation}}
\newcommand{\en}{\end{equation}}

\newcommand{\R}{\mathbb{R}}

\renewcommand{\vec}[1]{\boldsymbol{#1}}

\newcommand*{\pd}[2]{\mathchoice{\frac{\partial#1}{\partial#2}}
  {\partial#1/\partial#2}{\partial#1/\partial#2}
  {\partial#1/\partial#2}}
\newcommand*{\od}[2]{\mathchoice{\frac{d#1}{d#2}}
  {d#1/d#2}{d#1/d#2}{d#1/d#2}}

\begin{document}
\numberwithin{equation}{section}


\title{A Translation of the T. Levi-Civita paper:\\
  \emph{Interpretazione Gruppale\\
    degli Integrali di un Sistema Canonico}\\[5mm]
\small Rend.\ Acc. Lincei, s. 3$^\text{a}$, vol.\ VII, 2$^\text{o}$ sem.\ 1899,
pp.\ 235--238}


\author{G. Saccomandi$^{a}$,
  R.~Vitolo$^{b}$\thanks{This work has been partially supported by the Universit\`a degli Studi
  di Perugia, Universit\`a del Salento, GNFM of INdAM.}, \\
  {\it\small $^a$Dipartimento di Ingegneria Industriale,}
  \\[-6pt]
  {\it\small  Universit\`{a} degli Studi di Perugia, 06125 Perugia;}\\
  {\it\small $^b$Dipartimento di Matematica,}
  \\[-6pt]
  {\it\small Universit\`{a} del Salento, Via per Arnesano, 73100 Lecce.}}

\maketitle
\begin{abstract}
  In this paper we provide a translation of a paper by T. Levi-Civita,
  published in 1899, about the correspondence between symmetries and
  conservation laws for Hamilton's equations. We discuss the results of this
  paper and their relationship with the more general classical results by
  E. Noether.
\end{abstract}

\section{Introduction}


Tullio Levi-Civita (1873--1941) has been one of the most important
mathematicians of the last century. Levi-Civita is best known for his work on
the absolute differential calculus and because he was the founder with Gregorio
Ricci-Curbastro of the subject now known as tensor analysis. Levi-Civita
published nearly two hundred papers about classical mechanics, hydromechanics,
thermodynamics, elasticity, the strength of materials, astronomy,
electromagnetism, optics, relativity and quantum mechanics and pure
mathematics. He was also the author of several treatises, and in collaboration
with Ugo Amaldi he wrote an important reference book in classical mechanics:
\emph{Lezioni di Meccanica Razionale} \cite{LCA}. This is a three volumes book
first published in 1922 and then revised for the last time, after the death of
Levi-Civita, by Amaldi in 1949.

In reading \cite{LCA} it is impossible not to notice that this book does not
contains any mention to the Noether's 1914 theorem \cite{NT}. Indeed at page
$98$ of the wonderful book \cite{Yvette} we read:

\begin{quotation}
  The surveys of Vizgin [1972] and Kastrup [1987] and our own research have
  yielded a surprisingly small number of references to Noether before 1950. In
  particular, we are still astonished by the absence of citations dealing with
  invariance and its related mathematics in the corpus of the then standard but
  now classical textbooks on the variational calculus. There is nothing on
  invariance problems in the treatises of Tonelli [1921], \textbf{Levi-Civita
    and Amaldi} [1923], Bliss [1925], Forsyth [1927], Ames and Murnaghan [1929]
  who treat cyclic coordinates but not general invariance properties, nor
  later, in that of Elsgolc [1952], originally written in Russian and
  translated into English in 1961, cited by Gel'fand and Fomin [1961], nor in
  Pars [1962].
\end{quotation}

For someone who knows the Levi-Civita scientific production it is truly hard
to believe that the non-citation of Noether's theorem in \cite{LCA} was
possible.  This fact is also clear from a bird's eye reading of
\cite{LCA}. The authors of this wonderful treatise consider in great detail,
for a book in classical Mechanics, the theory of transformations, applying it
to different topics. They also refer to the fundamental papers of Sophus
Lie\footnote{To be more precise about this point, in the second part of the
  second volume of \cite{LCA} you can find a bio sketch of Sophus Lie at page
  302 and he is cited in 9 pages.}.  Moreover, both Levi-Civita and Amaldi
have been very active in studying the theory of transformation groups applied
to various problems of mathematical physics. For example in a recent issue of
the journal \emph{Regular and Chaotic Dynamics} Sergio Benenti has provided
the English translation of the Levi-Civita paper \emph{Sulle trasformazioni
  delle equazioni dinamiche} \cite{LC2}. It is interesting to report a comment
by the editorial board of \emph{Regular and Chaotic Dynamics} on \cite{LC2}:

\begin{quotation}
  This paper is truly a pioneering work in the sense that the real power of
  covariant differentiation techniques in solving a concrete and highly
  nontrivial problem from the theory of dynamical systems was demonstrated. The
  author skillfully operates and weaves together many of the most advanced (for
  that times) algebraic, geometric and analytic methods. Moreover, an attentive
  reader can also notice several forerunning ideas of the method of moving
  frames, which was developed a few decades later by E. Cartan.
\end{quotation}

We are not experts in the history of mathematics and we have not got a
complete and detailed knowledge of all the writings of Levi-Civita and Amaldi
(for example of all the letters that they have written and for which there is
still a record), but it seems to us that there are no documents of these two
mathematicians where the question of the theorem of Noether is discussed.  The
journal where \cite{NT} is published is cited in the famous text \cite{RLC01},
so it is plausible to think that Levi-Civita would read it more or less
regularly.

We looked at the collected papers of Levi-Civita in an effort to try to
understand possible reasons for this fact\footnote{Clearly, here we are
  approaching the problem in a purely speculative framework. It is possible
  that Levi-Civita and Amaldi have never read the paper by Noether and they
  were simply unaware of the result.}, and we found a very interesting short
note written at the onset of his career: the paper \cite{LC1}.

Our aim is to provide an English translation of \cite{LC1} because we think
that in this short note the connection between infinitesimal generators of
symmetries (in the sense of Sophus Lie) and first integrals has been proposed
for the first time in its general setting.

Let us stress in a clear and definitive way that we are not claiming that in
\cite{LC1} you can find the two Noether's theorems from \cite{NT}.  First of
all Levi-Civita considers Hamiltonian finite dimensional systems and the
connection between Lagrangian and Hamiltonian formulation was not rigorously
stated. Moreover, Noether's results are much more general than what is
contained in the short note by Levi-Civita.  What we are claiming is that the
connection between Lie's theory of groups and conserved quantities in
mechanics was already clear in \cite{LC1}.

To support our claim in Section~\ref{sec:description} we provide a paraphrase
of the paper \cite{LC1} using a modern language (see, \emph{e.g.},
\cite{BCA10,Many99,Olver}). We think that this \emph{exercise} is helpful to
clarify our point of view. We point out that in \cite{LC1} there are very few
references, among which we found a citation to two short notes by Maurice
L\'evy and a short note of Valentino Cerruti. In those references, integrals of
motion are in some sense connected with the symmetries of the mechanical
system, but in an episodic and incomplete way.

We end up this introduction with two remarks.

First of all we think that it is impossible to have a reasonable answer to our
starting question: why are Noether theorems not cited in \cite{LCA}? We
believe that the general problems about the inception and the reception of the
Noether theorems before 1950 stressed out in \cite{Yvette} in the case of
Levi-Civita have to be summed up to the following specific facts:
\begin{itemize}
\item The connection between symmetries and constant of motions was already
  clear to Levi-Civita. This means that maybe Levi-Civita was not so surprised
  by the result of Noether and therefore maybe he was not so interested in
  reading the paper.
\item Levi-Civita was interested more in the Hamiltonian formalism than in the
  Lagrangian one\footnote{Ugo Amaldi took the decision to stop his research
    activity soon in his career and therefore was no more up to date on the
    scientific literature to help Levi-Civita in discovering interesting
    papers.}.  We recall that Olver says in \cite{Olver} that: \emph{The
    Hamiltonian version of Noether's theorem has a particularly attractive
    geometrical flavour, which remains somewhat masked in our previous
    Lagrangian framework}.  We know that Levi-Civita was one of the first
  mathematicians interested in the geometry of Mechanics, and it is possible
  that he was more interested in this aspect that in the Lagrangean setting.
\end{itemize}

The second remark is about the fact that papers \cite{LC1}, \cite{LC2} and
many other papers by Levi--Civita seem to have not properly been noticed in
the literature.  Because of the racist laws introduced in Italy by the fascist
regime, the notice of the death of Levi--Civita was only given by the
Osservatore Romano and by foreign academies\footnote{Levi-Civita in 1938 was
  forced to retirement due to the promulgation of laws against Jews by the
  fascist regime.}. For a long time a complete list of the publications by
Levi--Civita was not available. Moreover, Ugo Amaldi's wish was that no
obituary of him had to be diffused after his death. Therefore, it is not
strange at all that the less famous papers of Levi--Civita have been forgotten
for a long time.

\section{Translation of the paper\protect\footnote{Presented by the Fellow
    \textsc{V. Cerruti} in the session of 5 November 1899}}
\label{sec:translation-paper}

Mr.\ \textsc{Maurice L\'evy} was the first to observe\footnote{\emph{Comptes
    Rendus}, t.\ LXXXVI, 18 February and 8 April 1878.} that a displacement
without deformation is possible in a generic manifold if and only if it is
possible to remove one of the variables from the square of the linear element
by a suitable transformation. This is equivalent to say that there exists an
homogeneous, linear first integral for the geodesics of the manifold.

Prof.\ \textsc{Cerruti} reexamined the topic\footnote{In this
  \emph{Rendiconti}, ser.\ 5$^\text{a}$, vol.\ III, 1895.} also considering
the case where conservative forces are acting. The above relation between
first integrals and rigid displacements can be formulated in the language of
groups as follows\footnote{See the Notes \emph{On the motion of a rigid body
    about a fixed point}, in this \emph{Rendiconti}, ser.\ 5$^\text{a}$, vol.\
  V, 1896 and the elegant proof of Mr. \textsc{Liebmann}, \emph{Math. Ann.},
  B. 50, 1897.}. If the kinetic energy and the potential admit the same
infinitesimal point transformation, then the equations of motion have an
homogeneous linear first integral, and vice versa. (The left-hand side of the
integral, as written in canonical form, coincides with the symbol of the
infinitesimal transformation).

A natural question arises: does any group-like character still correspond to
non-linear integrals?

The answer is affirmative and this holds for any canonical system
\begin{equation}\label{eq:S}\tag{S}
  \left\{\begin{array}{l}\displaystyle
\frac{dx_i}{dt}=\frac{\partial H}{\partial p_i} \\[5mm]
\displaystyle\frac{dp_i}{dt}=-\frac{\partial H}{\partial x_i}
\end{array}\right.
\qquad (i=1,2, \ldots, n),
\end{equation}
as long as one does not only consider point transformations (with respect to
$x$, operating on $p$ by prolongation), but more general contact
transformations of $x$, $p$. One finds indeed that \emph{integrals of a
  canonical system and contact transformations in $x$, $p$ changing the system
  into itself are essentially the same thing. For each integral there exists a
  transformation and conversely. The characteristic functions of the
  transformations \emph{(by suitably fixing an addend that remains apriori
    undetermined)} can be made to coincide with the left-hand sides of the
  corresponding integrals.}

The theorem is proved in a very simple way. Let
\begin{displaymath}
  \delta f = \xi_1\pd{f}{x_1}+\cdots+\xi_n\pd{f}{x_n}
  +\pi_1\pd{f}{p_1}+\cdots+\pi_n\pd{f}{p_n}
\end{displaymath}
be an infinitesimal transformation in $x$, $p$. Suppose that the increments
$\xi$, $\pi$ are functions of $x$, of $p$ and of a parameter $t$, which is
invariant under the transformation. By considering $x$, $p$ as functions of $t$
we can extend $\delta f$ to the derivatives $\od{x_i}{t}$, $\od{p_i}{t}$, and
the respective increments will be obtained by the formulae
\begin{gather*}
  \delta\od{x_i}{t}=\od{\delta x_i}{t}=\od{\xi_i}{t},\\
  \delta\od{p_i}{t}=\od{\delta p_i}{t}=\od{\pi_i}{t}.
\end{gather*}
Applying the transformation $\delta f$ to system~\eqref{eq:S} yields
\begin{equation}
  \label{eq:1LC}\tag{1}
    \left\{\begin{array}{l}\displaystyle
\delta\left\{\od{x_i}{x}-\pd{H}{p_i}\right\}=0 \\[5mm]
\displaystyle
\delta\left\{\od{p_i}{t}+\pd{H}{x_i}\right\}=0
\end{array}\right.
\qquad (i=1,2, \ldots, n).
\end{equation}
The above identities must be satisfied by virtue of \ref{eq:S} if
the system admits the infinitesimal transformation $\delta f$.

Let us introduce the hypothesis that $\delta f$ be a contact
transformation. The $\xi$ and the $\pi$ are derivatives of the same function
$W(x,p,t)$\footnote{\textsc{Lie--Engel}, \emph{Theorie der
    Transformationsgruppen}, vol.\ II, cap.\ 14.}, according to
\begin{equation}
  \label{eq:2LC}\tag{2}
  \xi_i=\pd{W}{p_i},\qquad \pi_i=-\pd{W}{x_i},
\end{equation}
and the symbol $\delta f$ becomes the Poisson bracket $(W,f)$.

Equation \eqref{eq:1LC} can be written as
\begin{gather*}
\od{\displaystyle\,\pd{W}{p_i}}{t} - \left(W,\pd{H}{p_i}\right)=0,\\
\od{\displaystyle\,\pd{W}{x_i}}{t} - \left(W,\pd{H}{x_i}\right)=0,
\end{gather*}
hence, evaluating derivatives and using \eqref{eq:S}:
\begin{gather*}
\pd{^2 W}{p_i\partial t}+\left(H,\pd{W}{p_i}\right)
  -\left(W,\pd{H}{p_i}\right)=0,\\
\pd{^2 W}{x_i\partial t}+\left(H,\pd{W}{x_i}\right)
  -\left(W,\pd{H}{x_i}\right)=0,
\end{gather*}
that, by virtue of bracket properties, is equivalent to
\begin{equation}
  \label{eq:3LC}\tag{1'}
  \left\{
  \begin{array}{l}
  \displaystyle\pd{}{p_i}\left[\pd{W}{t}+(H,W)\right]=0,\\[5mm]
  \displaystyle\pd{}{x_i}\left[\pd{W}{t}+(H,W)\right]=0.
  \end{array}
  \right.
\end{equation}
From the above equations it follows that $\pd{W}{t}+(H,W)$ depends on $t$
only. Now $W$, the characteristic function of $\delta f$, is determined
by~\eqref{eq:2LC} up to an additive function of $t$. One can always choose
this in such a way that the following identity holds:
\begin{equation}
  \label{eq:4LC}\tag{1"}
  \pd{W}{t}+(H,W)=0.
\end{equation}
Clearly, one can also obtain~\eqref{eq:1LC} starting from~\eqref{eq:4LC}
and tracing each step backwards.
Therefore~\eqref{eq:4LC} is a necessary and sufficient condition for the
canonical system~\ref{S} to admit the infinitesimal contact transformation
$(W,f)$.

On the other hand~\eqref{eq:4LC} precisely states that $W=\text{const}$ is
an integral of system~\eqref{S}. This proves the statement formulated above.

Note that if $W$ is linear and homogeneous in $p$ (and in this case only),
then $\delta f$ comes from the prolongation of a point transformation with
respect to $x$. It follows that the existence of a homogeneous linear integral
and the existence of a point transformation changing the canonical system into
itself are concomitant facts. In particular, the theorem by
\textsc{L\'evy--Cerruti} is reobtained by supposing that $H=T-U$ with T
homogeneous of second degree in $p$ and that $U$ is a function of $x$
only. It follows indeed from~\eqref{eq:4LC}, by splitting terms of different
degree in $p$, that $T$ and $U$ separately admit the transformation $W$.

\section{Some comments on the Levi-Civita paper}
\label{sec:description}

The paper by Levi-Civita considers an Hamiltonian system \be \label{S}
\frac{dx_i}{dt}=\frac{\partial H}{\partial p_i}, \quad
\frac{dp_i}{dt}=-\frac{\partial H}{\partial x_i}, \quad i=1, \ldots, n, \en
where $x_i=x_i(t)$ and $p_i=p_i(t)$. In a modern geometric language we can say
that the phase space is of the form $\R\times T^*M$, with coordinates
$(t,x_i,p_i)$, where $\R$ represents time and $M$ is an $n$-dimensional
manifold representing positions.  Now, it is clear that $F(t, x_i, p_i)$ is a
first integral for \eqref{S} if and only if
\[
\frac{\partial F}{\partial t}+ \{F,H\}=0
\]
on the manifold described in the phase space by the solutions of
\eqref{S}. Here, as usual, $\{\cdot,\cdot\}$ are the \emph{Poisson's brackets}
which if we are working in $\R^{2n}$ are defined as
\begin{equation}\label{eq:1}
\{F,H\}=
\frac{\partial F}{\partial x_i}\frac{\partial H}{\partial
    p_i}-\frac{\partial F}{\partial p_i}\frac{\partial H}{\partial x_i}
\end{equation}
(here and in what follows sums on repeated indexes are understood).  Then
Levi-Civita considers an infinitesimal transformation $\delta f$ for the
variables $x_1, \ldots,$ $x_n$ and $p_1, \ldots, p_n$ as dependent variables
with respect to the independent variable $t$. In other words, his
transformation preserves the projection $\R\times T^*M \to \R$. This is clear
when he speaks of the \emph{incrementi} of the variables $x$ and $p$, denoted
as $\xi_i$ e $\pi_i$. So, his infinitesimal transformation $\delta f$ can be
regarded as the time dependent vector field
\begin{equation}\label{eq:2}
  \vec{v}=\xi_i \frac{\partial}{\partial x_i}+\pi_i
\frac{\partial}{\partial p_i}
\end{equation}
on $\R\times T^*M$.  In this infinitesimal transformation Levi-Civita assumes
that time is invariable, \emph{i.e.} the coefficient of $\pd{}{t}$ is null.  In
other words, $\vec v$ is a \emph{vertical} vector field.

By introducing new coordinates $\dot x_i$ and $\dot p_i$ for
the derivatives $\od{x}{t}$ and $\od{p}{t}$ we see that the \emph{incrementi}
of those variables in Levi-Civita's language are just the coefficients of the
first prolongation of $\vec v$ to the first jet space of the projection
$\R\times T^*M \to \R$:
\begin{displaymath}
\text{pr}^{(1)} \vec{v}=\xi \frac{\partial}{\partial
    x_i}+\pi_i\frac{\partial}{\partial p_i}+\frac{d\xi_i}{dt}\pd{}{\dot
    x_i}+\frac{d\pi_i}{dt}\pd{}{\dot p_i}.
\end{displaymath}
In the above formula we are considering $d \cdot /dt $ as a total
derivative:
\begin{displaymath}
  \od{}{t}=\pd{}{t}+\dot x_i\pd{}{x_i}+\dot p_i\pd{}{x_i}
\end{displaymath}
(in \cite{Olver} the notation $D_t$ is used). The fact that the expression for
coefficient functions of $\text{pr}^{(1)} \vec{v}$ is very simple is related
with the fact that we are not transforming the independent variable $t$.

If we rewrite \eqref{S} as
\begin{displaymath}
E^{(1)}_i:=\frac{dx_i}{dt}-\frac{\partial H}{\partial p_i}=0, \quad
E^{(2)}_i:=\frac{dp_i}{dt}+\frac{\partial H}{\partial x_i}=0, \quad i=1,
\ldots, n
\end{displaymath}
the infinitesimal invariance criterion is given by
\begin{equation}
  \label{invarianza}
  \left\{\begin{array}{l}
\text{pr}^{(1)} \vec{v}(E^{(1)}_i)=0,
\\
\text{pr}^{(1)}
\vec{v}(E^{(2)}_i)=0,
\\
E^{(1)}_i=0,
\\
E^{(2)}_i=0,
\end{array}
\right.
\end{equation}
for $i=1,\ldots,n$.  Levi-Civita requires that $\vec v$ is a \emph{contact
  transformations} (see \cite{Lie} or \cite{Eis}, for example). More precisely,
$\vec v$ is a parameter-dependent infinitesimal contact transformation, the
parameter being $t$ as stated before. Clearly, this kind of transformation
should not be confused with contact transformations of $\R\times T^*M$ where
$t$ is the dependent variable, where $\pi_i$ is defined as the total derivative
of $W$ (see \cite{BCA10}).

Let us recall some basic facts on infinitesimal homogeneous contact
transformations. Infinitesimal homogeneous contact transformations of $T^*M$
are vector fields $X$ on $T^*M$. Their flow $F_s$ maps a cotangent
vector $(x_i,p_i)$ into another cotangent vector $(x'_i,p'_i)$ in such a way
that, if $p_i=(d_{(x_i)}f)_i$, where $f\colon M\to \R$, then
$p'_i=(d_{(x'_i)}f')_i$, where $f'=f\circ F_s\colon M\to \R$.
This characterization is equivalent to the requirement that
the flow preserves the contact form (or Liouville form) $\theta=p_idx_i$ on
$T^*M$, or that the equation $L_X\theta = 0$ (here $L_X$ stands for the Lie
derivative) holds. That implies the conditions
\begin{displaymath}
  \pi_j=-p_i\pd{\xi_i}{x_j},\qquad p_i\pd{\xi^i}{p_j}=0.
\end{displaymath}
It follows that, if $W=p_i\xi^i$, then $\xi_i=\pd{W}{p_i}$ and
$\pi_i=-\pd{W}{x_i}$. The function $W$ is uniquely defined.

Note that more general non-homogeneous contact
transformations can also be considered. A modern treatment of the subject can
be found in \cite[Chapter 2]{Many99}.

In our case we construct a function $W_t(x_i,p_i)$ for any value of the
parameter $t$, hence we obtain a time-dependent function $W(t,x_i,p_i)$
which is clearly defined up to an arbitrary function of time.

Using this fact the \eqref{invarianza} reads
\be \label{inv}
\frac{d}{dt}\left(\frac{\partial W}{\partial p_i} \right)-\left\{W,
  \frac{\partial H}{\partial p_i}\right\}=0, \quad
\frac{d}{dt}\left(\frac{\partial W}{\partial x_i} \right)-\left\{W,
  \frac{\partial H}{\partial x_i}\right\}=0.
\en
Consider the first equation
\eqref{inv} i.e
\[
\frac{\partial^2 W}{\partial p_i \partial t}+\frac{\partial^2 W}{\partial
  p_i \partial x_j} \frac{dx_j}{dt}+\frac{\partial^2 W}{\partial p_i \partial
  p_j} \frac{dp_j}{dt}-\left\{W, \frac{\partial H}{\partial p_i}\right\}=0,
\]
on the solutions of \eqref{S} we have
\[
\frac{\partial^2 W}{\partial p_i \partial t}+\left\{H, \frac{\partial
    W}{\partial p_i}\right\}-\left\{W, \frac{\partial H}{\partial
    p_i}\right\}=0.
\]
A similar computation for the second equation in \eqref{inv} says
\[
\frac{\partial^2 W}{\partial x_i \partial t}+\left\{H, \frac{\partial
    W}{\partial x_i}\right\}-\left\{W, \frac{\partial H}{\partial
    x_i}\right\}=0.
\]
We have shown
\be
\frac{\partial}{\partial p_i} \left[\frac{\partial
    W}{\partial t}+\left\{H, W\right\} \right]=0, \quad
\frac{\partial}{\partial x_i} \left[\frac{\partial W}{\partial t}+\left\{H,
    W\right\}\right]=0
\en
and being $W$ defined up to an arbitrary function of
time we have the desired result.

It is possible to recast the above discussion using a geometric framework, but
we think that our simple discussion is sufficient to clarify that in \cite{LC1}
we can find a clear and complete connection between symmetries of the Hamilton
equations (according with Lie theory) and first integral of finite dimensional
systems.

Of course, Noether's theorem \cite{NT} is formulated for Euler--Lagrange
equations in field theory. On the other hand, we could use a geometric
viewpoint to show that in mechanics, if the Lagrangian is regular\footnote{this
  means that the Legendre transformation is a local diffeomorphism}, then the
integral curves of the Euler--Lagrange equations and of the Hamilton equations
are (at least locally) in bijection (see, \emph{e.g.}, \cite[p.\
218]{AbrMar78}). Hence, the symmetries of the two equations are (at least
locally) in bijection. Of course, such a correspondence between symmetries
makes sense only if we allow symmetries to be generalized (or higher) (see,
\emph{e.g.}, \cite{Olver,Many99}), in such a way that we can consider
symmetries of the Euler--Lagrange equation that depend on velocities. On the
other hand, it is known from Noether's paper \cite{NT} generalized symmetries
allow for a complete identification between symmetries and conserved
quantities. So, generalized symmetries must be taken into account if we wish to
state Noether's theorem in the greatest generality possible.

Deeper discussions on the relationship between symmetries of the
Euler--La\-grange equations and the Hamilton equations with different
connections with Noether's theorem can be found in the papers
\cite{But06,GP92,MS77,RL82,SC81,SR02} in mechanics and \cite{GitTyu06} in field
theory.

Our starting idea was to understand why the theorem of Noether is not cited in
the treatise of rational mechanics by Levi-Civita.  Clearly we cannot provide a
definite answer to this question, but we have discovered the short note by
Levi-Civita that clearly adds interesting information to the history of the
correspondence among symmetries and conservation of laws.

\bigskip

\textbf{Acknowledgements.} We thank Yvette Kosmann-Schwarzbach, Giusep\-pe Gaeta
and Renato Vitolo for their interest in our work and for many useful comments
that helped us to improve the paper.



\end{document}